\documentclass[12pt]{article}
\usepackage{latexsym}
\usepackage{theorem}
\usepackage{graphicx}
\usepackage{amsmath}
\usepackage{color}
\usepackage{xcolor}
\usepackage{amsfonts}
\usepackage{natbib}
\usepackage{soul}
\usepackage{hyperref} 

\theorembodyfont{\rmfamily}
\newtheorem{theorem}{Theorem}

\newtheorem{lemma}[theorem]{Lemma}

\theoremstyle{break}

\newcommand{\la}{\lambda}

\newcommand{\new}[1]{{\color{black}#1}}

\graphicspath{{../fig/}}

\title{Optimal Patrol of a Perimeter}

\author{Kyle Y. Lin\thanks{Operations Research Department, Naval Postgraduate School, Monterey, CA 93943, kylin@nps.edu}}

\date{September 11, 2020}

\providecommand{\keywords}[1]{\textbf{\textbf{Keywords:}} #1}

\begin{document}

\maketitle

\begin{abstract}
\noindent
A defender dispatches patrollers to circumambulate a perimeter to guard against potential attacks. The defender decides on the time points to dispatch patrollers and each patroller's direction and speed, as long as the long-run rate patrollers are dispatched is capped at some constant. An attack at any point on the perimeter requires the same amount of time, during which it will be detected by each passing patroller independently with the same probability. The defender wants to maximize the probability of detecting an attack before it completes, while the attacker wants to minimize it. We study two scenarios, depending on whether the patrollers are undercover or wear a uniform. Conventional wisdom would suggest that the attacker gains advantage if he can see the patrollers going by so as to time his attack, but we show that the defender can achieve the same optimal detection probability by carefully spreading out the patrollers probabilistically against a learning attacker.
\end{abstract}


\section{Introduction}
\label{sec:introduction}
\new{
Patrol is a common practice to secure the perimeter of an area of interest, such as soldiers patrolling the perimeter of a military installation to guard against illicit activities.
With limited manpower, how should the patrol schedule be set in order to maximize patrol effectiveness?
If the patrollers wear a uniform, can an attacker gain advantage by observing the patrol process so as to time his attack?
This paper seeks to answer these questions with a game-theoretic model.
}

Consider a two-person zero-sum Stackelberg game played by a defender and an attacker on an infinite time horizon $[0, \infty)$.
The defender moves first and is tasked to secure the perimeter of an area of interest by dispatching patrollers from a base located on the perimeter.
Each patroller departs from the base, circumambulates the perimeter either clockwise or counterclockwise at possibly variable speed, and returns to the base to complete a full circle.
The defender decides on the patrol schedule---the time points at which patrollers depart from the base, the moving direction of each patroller, and each patroller's speed---as long as the long-run rate at which patrollers depart from the base \new{is equal to} $\la > 0$.
Specifically, if an increasing sequence $s_1, s_2, \ldots$ represents the time points patrollers are dispatched, then it is a feasible pure strategy if
\new{
\begin{align}
\lim_{T \rightarrow \infty} \frac{\max \{n \in \mathbb{N} \, | \, s_n \leq T\}}{T} = \la,
\label{eq:lambda}
\end{align}
where $\mathbb{N}$ is the set of positive integers and $\max \{n \in \mathbb{N} \, | \, s_n \leq T\}$ is the number of patrollers dispatched by time $T$.}
Any probability measure on the set of pure strategies is a feasible mixed strategy for the defender.
For example, the defender can dispatch patrollers according to a Poisson process at rate $\la$, or at fixed intervals of $1/ \la$ time units, and ask each patroller to go in the same clockwise direction at the same constant speed.

The attacker moves second and seeks to carry out an attack on the perimeter---such as penetrating the perimeter or planting a hidden bomb.
\new{The attacker chooses which point on the perimeter to attack and when to show up, and can observe the patrol process for as long as he wants before starting an attack.
Specifically, the attacker is free to choose \textit{any} operating window $[y, z]$, as long as $0 < y < z$.}
Upon arrival at time $y$, the attacker starts to observe the patrollers going by and decide when to initiate an attack, as long as the attack completes by time $z$.
An attack takes a fixed amount of time $t >0$, during which the attack will be detected by each passing patroller independently with probability $p \in (0, 1]$.
The attacker's objective is to minimize the probability of getting detected by any patroller during his attack, while the defender's objective is the opposite.
We denote this \new{two-person zero-sum Stackelberg} game by $G(\la, t, p)$.

The patrol game $G(\la, t, p)$ captures the essence of a real-world situation where patrollers can be easily recognized by the attacker, such as ground robots, \new{police cars}, and soldiers in uniform.
In some other situations, patrollers may be invisible to the attacker, such as undercover patrol agents.
Conventional wisdom would suggest that, by seeing the patrollers, the attacker can gain advantage by learning about the patrol pattern in real time and use that knowledge to improve his chance.
Somewhat surprisingly, we show that the defender can achieve the same detection probability---whether or not the patrollers are visible---by carefully spreading out the patrollers across time in a probabilistic manner.

Patrol problems arise in many real-life situations.
Police officers patrol cities and highways; security guards patrol museums and airport lobbies; soldiers patrol military bases and borders.
Earlier works focus on allocating police patrol resources among different areas to maximize the overall performance~\citep{Chaiken:1978p5794,Chelst:1978p5733,Larson1972,Olson:1975p5791,szechtman2008models,Portugal20131572}.
Besides police patrol in urban areas, there are specialized patrol models for rural areas \citep{Birge:1989p5790} and on highways \citep{Lee:1979p5793,Taylor:1985p5797}.
These earlier works assume that the frequencies of crimes at different locations are \new{exogenous and} known to the patrol force.



In recent years, there has been a growing interest in taking a game-theoretic approach to study patrol problems.
\cite{Alpern:2011} study a patrolling game played on a graph, where one patroller traverses the graph via its edges to detect possible attacks at vertices, while an attacker chooses which vertex to attack to minimize the probability of getting detected.
\cite{Lin:2013} extends the work to allow random attack times, and \cite{Lin2014nrl} further accounts for the possibility of overlook.
\cite{papadaki2016patrolling} and \cite{alpern2019optimizing} study a similar formulation but focus on the border patrol, while \cite{mcgrath2017robust} consider how to coordinate multiple patrollers.
In these works, the attacker cannot see the patroller.


Patrollers that are visible to the attacker when near are often considered in the context of robot patrol.
\cite{Basilico:2009uq,basilico2012patrolling} study a graph patrolling game in which the intruder  sees the patroller and uses that knowledge to decide when and where to initiate an attack.
\cite{NoaAgmon:2008td,agmon2008impact} place multiple robots on a perimeter and randomize their movements so that the intruder cannot predict the robots' locations with certainty.
\cite{zoroa2012patrolling} also study perimeter patrol, but in their model, the intruder can move along the perimeter during an attack while the patroller can choose any point on the perimeter to inspect in each time period.
\new{\cite{alpern2019note} considers a uniformed patrolling game on a star network between one patroller and one attacker, and shows that the attacker can improve the probability of a successful attack by timing his attack based on when the patroller was last seen.}


The vast majority of earlier patrol models discretize time and space for mathematical tractability.
In this paper, we treat time and space as continuum to improve model realism.
This approach is adopted in \cite{alpern2016patrolling} and \cite{Garrec:2019}, which study patrol problems in a different context with the patrollers invisible to the attacker.
\cite{Ruckle:1983p5984} studies geometric games, which shares similar spirit to our patrol game, but with different strategy sets and payoff functions.

\new{The rest of this paper proceeds as follows.
Section~\ref{sec:location} focuses on the interaction between the attacker and patrollers at the location where the attacker chooses to attack.
Section~\ref{sec:perimeter} presents the solution to the perimeter patrol game $G(\la, t, p)$.
Section~\ref{sec:finite} discusses a patrol game played on a finite time horizon and its connection with $G(\la, t, p)$.
Section~\ref{sec:concluding} offers a few concluding remarks.
}


\section{Patrol Game at One Location}
\label{sec:location}
\new{In order to analyze the perimeter patrol game $G(\la, t, p)$, we first focus on the interaction between the attacker and the patrollers at the location on the perimeter where the attacker chooses to attack.
In addition, we will restrict the defender to send all patrollers in the same direction at the same constant speed, so a pure strategy for the defender can be reduced to the time points at which patrollers arrive at the location.
Section~\ref{sec:visible} formulates and solves this patrol game played at one location with visible patrollers, and Section~\ref{sec:invisible} shows that the value of the game remains the same even if the patrollers are invisible to the attacker.}

\subsection{Visible Patrollers}
\label{sec:visible}
Consider a modification to the game $G(\la, t, p)$ by focusing our attention at a single location on the perimeter.
\new{
The defender moves first and chooses time points $s_1, s_2, \ldots$, at which patrollers arrive at the location.
Any increasing sequence $s_1, s_2, \ldots$ that satisfies \eqref{eq:lambda} is a feasible pure strategy for the defender.
The attacker moves second and is free to choose any operating window $[y, z]$, as long as $0 < y < z$.
Starting at time $y$, the attacker observes the patrollers going by and decides when to initiate an attack, but has to complete the attack by time $z$.
An attack takes a fixed amount of time $t >0$, during which the attack will be detected by each arriving patroller independently with probability $p \in (0, 1]$.
The attacker's objective is to minimize the probability of getting detected by any patroller during his attack, while the defender's objective is the opposite.
Denote this two-person zero-sum Stackelberg game by $\Gamma(\la, t, p)$.
}

\new{If the expected number of the arriving patrollers during an attack is fixed, how does the defender randomize the number to maximize the detection probability?
This questions is key to solving $\Gamma(\la, t, p)$ and is answered in Lemma~\ref{le:N}, whose proof is reproduced from \cite{Lin2014nrl}.
}

\begin{lemma}
\label{le:N}
Denote the number of arriving patrollers during an attack by $N$, a nonnegative integer-valued random variable, with each patroller detecting the attack independently with probability $p \in (0,1]$.
If $E[N] =c >0$ is fixed constant, then the distribution that maximizes the detection probability $1 - E[ (1-p)^N]$ is the one that has the least variance.
In other words, if $c$ is an integer, then it is optimal to take $N=c$ deterministically; otherwise, the optimal distribution is
\begin{align*}
P\{N= \lfloor c \rfloor\} = \lceil c \rceil - c, \\
P\{N= \lceil c \rceil\} = c - \lfloor c \rfloor.
\end{align*}
\end{lemma}
\textit{Proof.}
Condition on $N$ to compute the objective function
\begin{align*}
E[1 - (1-p)^N] &= \sum_{n=0}^\infty  (1 - (1-p)^n) P\{N=n\}  \nonumber \\
&= \sum_{n=1}^\infty \left( p \sum_{k=0}^{n-1} (1-p)^k P\{N=n\} \right) \nonumber \\
&= p \sum_{k=0}^\infty \sum_{n=k+1}^\infty (1-p)^k P\{N=n\} \nonumber \\
&= p \sum_{k=0}^\infty (1-p)^k P\{N>k\}.
\end{align*}

Recall that $\sum_{k=0}^\infty P\{N>k\} = E[N] = c$ is fixed constant.
Writing $x_k = P\{N>k\}$ as decision variables, $k=0,1,\ldots$, the optimization problem can be rewritten as
\begin{align*}
\max & \quad \sum_{k=0}^\infty (1-p)^k x_k, \\
\text{subject to} & \quad \sum_{k=0}^\infty x_k = c, \\
& \quad 1 \geq x_0 \geq x_1 \geq \cdots \geq 0.
\end{align*}
Because $(1-p)^k$ decreases in $k$, the optimal solution to the preceding problem is to let
\[
x_k = \left\{
\begin{array}{ll}
1, & k=0,1,\ldots, \lfloor c \rfloor -1, \\
c - \lfloor c \rfloor, & k = \lfloor c \rfloor, \\
0, & k \geq \lfloor c \rfloor + 1.
\end{array}
\right.
\]
In other words, the optimal choice is for $N$ to take on the two integers surrounding $E[N]$, or just $E[N]$ if it happens to be an integer.
\hfill
$\Box$

\bigskip

\new{Lemma~\ref{le:N} can also be proved by observing that the detection probability $f(n) = 1 - (1-p)^n$ is a concave function in the number of arriving patrollers $n$, for $n=0,1,\ldots$, so $E[f(N)]$ is maximized by minimizing the variance of $N$ when $E[N]$ is fixed.}
Intuitively, because the expected number of patrollers encountered during an attack period is fixed, and detecting the attack twice or more times is just as good as detecting the attack once, it is better for each attack to encounter roughly the same number of patrollers rather than some encountering many and some very few.

\new{In accordance with the spirit of Lemma~\ref{le:N}, a reasonable strategy for the defender in the patrol game $\Gamma(\la, t, p)$ is to send patrollers at fixed intervals of $1/\la$ time units, so that an attack will encounter either $\lfloor \la t \rfloor$ or $\lceil \la t \rceil$ patrollers, regardless of when the attack begins.
The downside of this strategy, however, is that the attacker can see the patrollers going by, so he can begin an attack immediately after a patroller has just left.
Such an attack will encounter only $\lfloor \la t \rfloor$ patrollers, so the resulting detection probability falls short of the upper bound in Lemma~\ref{le:N}.
It turns out, by properly randomizing the patrol pattern, it is possible to achieve the desired distribution in Lemma~\ref{le:N}, as seen in the next theorem.}

\begin{figure}[t]
\begin{center}
\includegraphics[scale=0.75]{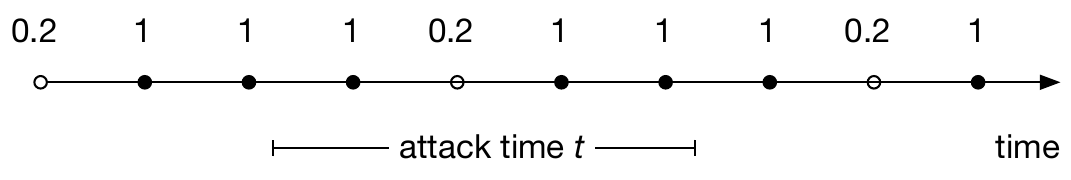}
\caption{The optimal patrol strategy presented in Theorem~\ref{th:main}, with $\la t=3.2$, so \new{$\lceil \la t \rceil = 4$} and $\Delta = t/4$.
Each solid dot corresponds to a \textit{blue} patroller, and each circle corresponds to a potential \textit{red} patroller with probability 0.2.
Regardless of when the attack begins, the attack period $t$ will encounter 3 patrollers with probability 0.8, or 4 patrollers with probability 0.2.
}
\label{fig:OptPatrol}
\end{center}
\end{figure}

\begin{theorem}
\label{th:main}
\new{The defender's optimal strategy in the Stackelberg game $\Gamma(\la, t, p)$ is as follows.
If $\la t$ is an integer, then place patrollers at fixed intervals of $1/\la$ time units.
If $\la t$ is not an integer, write $\Delta = t / \lceil \la t \rceil$.
Place a patroller at time points $j t + k \Delta$ for $j=0, 1, 2, \ldots$ and $k=1, 2, \ldots, \lfloor \la t \rfloor$; and a patroller at time points $j t$ independently with probability $\la t - \lfloor \la t \rfloor$, for $j=0, 1, 2, \ldots$.
The value of the game is given by
\begin{equation}
V(\la, t, p) = 1 - (\la t - \lfloor \la t \rfloor) (1-p)^{\lceil \la t \rceil} - (1 - \la t + \lfloor \la t \rfloor) (1-p)^{\lfloor \la t \rfloor},
\label{eq:value}
\end{equation}
which increases in $\la$, $t$, and $p$.}
\end{theorem}
\textit{Proof.}
\new{We first show that the proposed strategy guarantees the detection probability in \eqref{eq:value}.
If $\la t$ is an integer, then regardless of when the attack begins, the attack period will encounter exactly $\la t$ patrollers, so the defender guarantees detection probability $1 - (1-p)^{\la t}$, which is equal to \eqref{eq:value} by taking $\lfloor \la t \rfloor = \lceil \la t \rceil = \la t$.
If $\la t$ is not an integer, call it a \textit{blue} patroller if the patroller arrives at time points $j t + k \Delta$ for $j=0, 1, 2, \ldots$ and $k=1, 2, \ldots, \lfloor \la t \rfloor$, and a \textit{red} patroller if the patroller arrives at time points $j t$ for $j=0,1,2,\ldots$.
An example with $\la t = 3.2$ is shown in Figure~\ref{fig:OptPatrol}.
Regardless of when the attack begins, any attack period will encounter exactly $\lfloor \la t \rfloor$ blue patrollers, and 1 red patroller with probability $\la t - \lfloor \la t \rfloor$.
In other words, each attack will encounter $\lceil \la t \rceil$ patrollers with probability $\la t - \lfloor \la t \rfloor$, or $\lfloor \la t \rfloor$ patrollers with probability $1-(\la t - \lfloor \la t \rfloor)$, so the defender guarantees detection probability in \eqref{eq:value}.}

We next show that the defender cannot increase the guaranteed detection probability in \eqref{eq:value}.
Recall from Lemma~\ref{le:N} that the detection probability in \eqref{eq:value} is the maximum the defender can achieve if $E[N] = \la t$, where $N$ is the number of arriving patrollers during the attack period $t$.
If the defender were to increase the guaranteed detection probability, she needs to increase $E[N]$, which is impossible because the attacker is allowed to choose the operating window $[y, z]$ arbitrarily.

\new{Finally, the monotonicity of $V(\la, t, p)$ in each of its three arguments can be verified with straightforward algebra.}
\hfill
$\Box$

\bigskip

\new{While the defender has a strategy to guarantee detection probability at least $V(\la, t, p)$ in \eqref{eq:value}, as shown in Theorem~\ref{th:main}, it is worth noting that the attacker does \textit{not} have a strategy to guarantee the detection probability no more than $V(\la, t, p)$, as explained below.
Recall that an attacker's pure strategy involves an operating window $[y, z]$.
Consider an arbitrary mixed strategy for the attacker, and write random variable $Z$ for the upper bound of the operating window induced by the mixed strategy.
For any $\epsilon > 0$, the defender can find some large number $b$, such that $P\{Z \leq b\} > 1 - \epsilon$; that is, the probability the attack will complete before time $b$ is at least $1 - \epsilon$.
}
Hence, the defender can jam up the time interval $[0, b]$ with patrollers to achieve detection probability $1 - \epsilon$, while stopping the patrol process after time $b$ altogether to still have the long-run rate capped by $\la$.
\new{Because for any attacker's strategy, the defender can find a patrol strategy to achieve detection probability arbitrarily close to 1, it follows that the attacker does not have a strategy to guarantee the detection probability no more than $V(\la, t, p)$.
Consequently, if the patrol game $\Gamma(\la, t, p)$ were formulated as a simultaneous-move game, then the game would not have a solution.
}


One may conjecture that dispatching patrollers according to a Poisson process would work well, because the attacker cannot learn about the future of the Poisson process from its past.
The downside of dispatching patrollers according to a Poisson process, however, is that some attacks would encounter many patrollers while some encounter very few or even none at all, which overall results in a lower detection probability than that in \eqref{eq:value}---according to Lemma~\ref{le:N}.
If the attack period is relatively short or the patrol resource is very limited, so that $\la t < 1$, then according to Theorem~\ref{th:main}, it is optimal for the defender to dispatch one patroller with probability $\la t$ after each $t$ time units, so that every attack---regardless of when it starts---will encounter one patroller with probability $\la t$ or no patroller with probaiblity $1 - \la t$.
The time between two consecutive patrollers is distributed as a geometric random variable with parameter $\la t$ multiplied by $t$.
In the limit as $t \rightarrow 0$, this interarrival time converges to the exponential distribution with rate $\la$, so the optimal patrol process converges to the Poisson process with rate $\la$.

\subsection{Invisible Patrollers}
\label{sec:invisible}
An interesting and important question is whether the defender can increase the detection probability by making the patrollers invisible to the attacker, or conversely, how much advantage the attacker gains by observing the patrollers going by to time his attack.
Somewhat surprisingly, the value of the game remains the same, \new{as seen in the next theorem.}


\begin{theorem}
Consider the patrol game $\Gamma(\la, t, p)$ but the patrollers are invisible to the attacker.
The defender's optimal strategy in Theorem~\ref{th:main} is still optimal, and the value of the game remains the same as in \eqref{eq:value}.
\end{theorem}
\new{
\textit{Proof.}
The defender can use the same strategy presented in Theorem~\ref{th:main} to guarantee the same detection probability in \eqref{eq:value}, but can she improve it?
Recall that the attacker is free to choose any operating window $[y, z]$.
Because the patrollers are invisible, a pure strategy for the attacker reduces to the time point, say $y$, at which  he begins the attack.
According to Lemma~\ref{le:N}, the detection probability in \eqref{eq:value} is the maximum the defender can achieve if $E[N] = \la t$, where $N$ is the number of arriving patrollers during the attack period $t$.
If the defender were to increase the guaranteed detection probability, she would need to increase $E[N]$, which is impossible because the attacker is free to choose any $y \geq 0$.
\hfill
$\Box$
}

\bigskip

In the case where the attacker cannot see the patrollers, the defender has a simpler patrol strategy that is also optimal:
Dispatch the first patroller at a random time according to the uniform distribution over $[0, 1/\la]$, and subsequent patrollers at fixed intervals of $1/\la$ time units.
Regardless of when the attack begins, the time elapsed between the start of the attack and the arrival time of the first patroller follows a uniform distribution over  $[0, 1/\la]$, so the detection probability remains the same as in \eqref{eq:value}.
This simpler patrol strategy, however, is not optimal if the attacker can see the patrollers, \new{unless $\la t$ is an integer.
If $\la t$ is not an integer and the patrollers arrive at fixed intervals of $1/\la$ time units, then the attacker can initiate an attack immediately after a patroller goes by so that he will encounter exactly $\lfloor \la t \rfloor$ but never $\lceil \la t \rceil$ patrollers.
To achieve the optimal detection probability, the defender must randomize the patrol schedule as described in Theorem~\ref{th:main}.
}




\section{Patrol Game on a Perimeter}
\label{sec:perimeter}
We now return to the patrol game on a perimeter $G(\la, t, p)$ introduced in Section~\ref{sec:introduction}.
To avoid triviality, the attacker cannot begin his operation until the first patroller has time to circumambulate the perimeter.
If we impose a constraint that each patroller circumambulates the perimeter in the same direction at the same constant speed before returning to the base, then it is inconsequential which point on the perimeter the attacker chooses to attack, because the interarrival times between two consecutive patrollers at any point on the perimeter is dictated by the difference between their departure times from the base.
Furthermore, the interaction between the patrollers and the attacker at any point on the perimeter is captured by the game $\Gamma(\la, t, p)$, so the value of the game remains the same as $V(\la, t, p)$ in \eqref{eq:value}.

Can the defender increase the detection probability if some patrollers circumambulate the perimeter clockwise and some others counterclockwise, or if they are allowed to move at different speeds or even change speed during a trip?
The answer is no.
As long as each patroller is not allowed to turn around during a patrol trip, each patroller goes by each point on the perimeter exactly once during a trip.
Hence, the long-run patrol rate observed at any point on the perimeter is still capped at $\la$, so the same analysis in Theorem~\ref{th:main} applies.
In other words, the optimal patrol strategy requires each patroller to move in the same direction at the same constant speed.
The choice of the patroller's (same constant) speed is inconsequential, because the time between consecutive patrollers passing by any point on the perimeter is dictated by the time between their departures from the base.
Hence, we have proved the following.

\begin{theorem}
Consider the patrol game on the perimeter $G(\la, t, p)$ described in Section~\ref{sec:introduction}.
The value of the game is given by \eqref{eq:value}.
The optimal strategy for the defender is to use the mixed strategy described in Theorem~\ref{th:main} to dispatch the patrollers from the base, and require each patroller to circumambulate the perimeter in the same direction at the same constant speed.
\end{theorem}

\new{In practice, the single-encounter detection probability $p$ may decrease if the patroller's speed increases.
In addition, due to manpower constraint, it may be necessary to require each patroller to circumambulate the perimeter under a certain amount of time.
In this situation, the defender should instruct each patroller to travel at the slowest feasible speed that meets the round-trip time constraint so as to make $p$ as large as possible.
}

\section{A Finite Time Game}
\label{sec:finite}
If the patrol game is played on a finite time horizon, and the defender is given a fixed number of patrollers, then by counting patrollers going by the attacker may be able to improve his probability of completing a successful attack.
In practice, the defender's planning horizon is usually much longer than that of the attacker.
For example, the defender may have 40 patrollers to dispatch during a 12-hour work day, while an attack takes 10 minutes to complete and the attacker can hide for up to 1 hour or so without raising suspicion.
This section describes a finite time game for this situation, and shows that the advantage gained by the attacker from observing the patrol process is rather limited.

Consider a modification of $\Gamma(\la, t, p)$, which is now played on a finite time horizon $[0, n^2 t]$, for some positive integer $n$.
The defender has exactly $k \equiv \lfloor \la  n^2 t \rfloor$ patrollers and chooses $k$ time points in $[0, n^2 t]$ to place them.
The attacker chooses an operating window $[y, y + n t]$; in other words, the attacker arrives at the location at time $y$ to observe the patrol process to time his attack, which must be completed by time $y + nt$.
The two players move \textit{simultaneously}, with the defender trying to maximize the detection probability while the attacker trying to minimize it.
We denote this \new{two-person zero-sum simultaneous-move game} by $\Gamma_n(\la, t, p)$ and bound its value in the next theorem.


\begin{theorem}
\label{th:finite}
Consider the finite time game $\Gamma_n(\la, t, p)$.
Write $m \equiv \lfloor k/n^2 \rfloor$ and $r \equiv k/n^2 - \lfloor k/n^2 \rfloor$, and note that $k/n^2 = m+r$ is the average number of patrollers in an interval of length $t$.
If $r = 0$, then the value of the game $V_n(\la, t, p)$ is
\begin{equation}
V_n(\la, t, p) = 1 - (1-p)^m.
\label{eq:r=0}
\end{equation}
If $r > 0$, then the value of the game satisfies the bounds
\begin{equation}
1 - q (1-p)^{m+1} - (1-q) (1-p)^m \leq V_n(\la, t, p) \leq 1 - r (1-p)^{m+1} - (1-r) (1-p)^m,
\label{eq:r>0}
\end{equation}
where 
\[
q \equiv \frac{n^2 r - (n-1)}{n^2 - (n-1)} < r.
\]
\end{theorem}
\textit{Proof.}
Consider an attack strategy that divides $[0, n^2 t]$ into $n^2$ equal-length subintervals of length $t$, and chooses each subinterval as the attack period uniformly randomly.
This strategy is feasible by choosing an operating window that covers the intended attack period and ignoring the patrol process observed, if any.
Write $N$ for the number of patrollers during the attack period.
Because the defender place $k$ patrollers and the attacker chooses each subinterval uniformly randomly with probability $1/ n^2$, we have that $E[N] = k / n^2$.
The upper bounds for $V_n(\la, t,p)$ in \eqref{eq:r=0} and \eqref{eq:r>0} follow immediately from Lemma~\ref{le:N}.

To find a lower bound for $V_n(\la, t, p)$, consider the following patrol strategy.
\begin{enumerate}
\item
$k/n^2$ is an integer, so $m=k/n^2 \geq 1$ and $r=0$.

Define $\Delta \equiv t/m$, so that the time horizon $[0, n^2 t]$ contains $k = n^2 m$ subintervals of length $\Delta$.
Generate a number $u$ from the uniform distribution over $(0, \Delta)$, and dispatch the $k$ patrollers at time points $u + i \Delta$, for $i=0, 1, \ldots, k-1$.
Because the consecutive patrollers are exactly $\Delta = t/m$ time units apart, regardless of the attacker's choice of attack period, the attack will encounter exactly $m$ patrollers.
Hence, this patrol strategy guarantees $V_n(\la, t, p) \geq 1 - (1-p)^m$, which proves \eqref{eq:r=0}. 
 
\item
$k/n^2$ is not an integer, so $r>0$.

Define $\Delta \equiv t / (m+1)$, so that the time horizon $[0, n^2 t]$ contains $n^2 t / \Delta = n^2 (m+1)$ subintervals of length $\Delta$.
Generate a number $u$ from the uniform distribution over $(0, \Delta)$, and dispatch $n^2 m$ patrollers at time points $u + i t + j \Delta$ for $i=0, 1, 2, \ldots, n^2-1$, and $j=1, 2, \ldots, m$.
Designate these $n^2 m$ patrollers as \textit{blue} patrollers and the other $k - n^2 m = n^2 r$ patrollers as \textit{red} patrollers.
Among the $n^2$ time points $u + i t$, for $i=0, 1, \ldots, n^2-1$, choose a subset of $n^2 r$ time points uniformly randomly to dispatch the $n^2 r$ red patrollers.

Recall that the attacker first chooses an operating window of length $n t = n (m+1) \Delta$ and can observe the patrol process during that window to time the attack.
Divide this operating window into $n$ equal-length subinterval, each having length $t$, we see that each subinterval contains exactly $m$ blue patrollers and possibly one red patroller.
The best case scenario for the attacker is to see a red patroller in every one of the first $n-1$ subintervals in the operating window, so that the probability that the last subinterval also contains a red patroller would be the smallest, namely,
\[
q = \frac{n^2 r - (n-1)}{n^2 - (n-1)}.
\]
Therefore, this patrol strategy guarantees the detection probability at least
\[
1 - q (1-p)^{m+1} - (1-q) (1-p)^m,
\]
which is the lower bound in \eqref{eq:r>0}, thus completing the proof.
\hfill $\Box$
\end{enumerate}

\bigskip


Attach a subscript $n$ to the variables used in Theorem~\ref{th:finite} to signify their dependency on $n$, so that $k_n \equiv \lfloor \la (n^2 t) \rfloor$, $m_n \equiv \lfloor k_n/n^2 \rfloor$, $r_n \equiv k_n/n^2 - \lfloor k_n/n^2 \rfloor$, and
\[
q_n \equiv \frac{n^2 r_n - (n-1)}{n^2 - (n-1)}.
\]
It is straightforward to verify that
\[
\lim_{n \rightarrow \infty} m_n = \lfloor \la t \rfloor, \qquad  \lim_{n \rightarrow \infty} q_n = \lim_{n \rightarrow \infty} r_n = \la t - \lfloor \la t \rfloor,
\]
and consequently, $\lim_{n \rightarrow \infty} V_n(\la, t, p) = V(\la, t, p)$ in \eqref{eq:value}, which is the value of the Stackelberg game $\Gamma(\la, t, p)$ in Theorem~\ref{th:main}.

While it is possible for the attacker to improve his chance if the defender has a fixed number of patrollers to dispatch in a finite time horizon, this advantage is limited if the defender's planning horizon is orders of magnitude higher than the attacker's operating window.
As $n$ increases, the dependency between red patrollers showing up at different time points diminishes, so does the attacker's ability to time his attack by observing the patrol process.

\section{Concluding Remarks}
\label{sec:concluding}
While we present the problem in the context of perimeter patrol, the model applies to any patrol problem in which the patrollers depart from a base, follow the same route, and return to the same base.
For example, police cars may patrol the same route in a city, leaving and returning to the same police station.
While the defender's optimal strategy may not be implemented to perfection---police officers cannot move at constant speed in traffic---the insight derived from the analysis nevertheless provides useful and important guideline on how to execute patrol scheduling.

If each passing patroller can detect an attacker in hiding with some probability, then it only makes the attacker's hide-and-wait strategy worse, so our results still hold.
If the attacker is allowed to wait, withdraw, and try again at a later time, the results still hold.
The optimal patrol strategy neutralizes any tactics the attacker may employ and produces the same detection probability regardless of when the attack begins.

Our model assumes that the detection probability $p$ and attack period $t$ remain the same for the entire perimeter, so the analysis offers great insight into real-world situations where these two quantities have minimal variation on the perimeter.
For example, police cars patrol night streets to bust drug dealing and prostitution, or unmanned vehicles use heat image to detect intrusions of a military installation.
If the attack period $t$ has a few different values---say 5 minutes to climb over a higher wall or 3 minutes over a lower wall---then it is sensible to set $t$ to be the smallest of these values to calculate the guaranteed detection probability.
If attacks are more difficult to see at certain locations, then one possibility is to instruct every patroller to slow down at those locations and to speed up at other locations---so as to keep the total time needed to circumambulate the perimeter under some specified amount---to bring $p$ to the same level for the entire perimeter.
A patrol problem with huge variation in $p$ and $t$ at many different points on the perimeter typically requires a different analysis that focuses on how to allocate more patrol resource to the weak spots, which is beyond the scope of this work.
Another possible future research direction is how to counter coordinated attack, if the attacker recruits an accomplice to observe or disrupt patrol patterns at another point on the perimeter.

\section*{Acknowledgments}
The author thanks anonymous referees and an associate editor for many constructive comments that substantially improve the mathematical rigor and exposition of this paper.

\bibliographystyle{apalike}
\bibliography{../bib/patrol}

\end{document}